\documentclass[10pt]{amsart}

\usepackage[unicode]{hyperref}
\usepackage{color}

\usepackage{amsmath,amsfonts,amssymb,amsthm}

\usepackage{longtable}

\numberwithin{equation}{section}

\newtheorem*{theorem*}{Theorem}

\newtheorem*{proposition*}{Proposition}

\newtheorem*{lemma*}{Lemma}

\newtheorem*{corollary*}{Corollary}

\theoremstyle{definition}

\theoremstyle{remark}

\makeatletter
\newcommand{\symbitem}[1]{\item[#1]%
\renewcommand{\@currentlabel}{#1}\ignorespaces}
\makeatother
\newcommand{\beq}{\begin{equation}}
\newcommand{\eeq}{\end{equation}}
\newcommand{\beqa}{\begin{eqnarray}}
\newcommand{\eeqa}{\end{eqnarray}}
\newcommand{\beaa}{\begin{eqnarray*}}
\newcommand{\ben}{\begin{eqnarray*}}
\newcommand{\eaa}{\end{eqnarray*}}
\newcommand{\een}{\end{eqnarray*}}

\def \D {\mathcal{D}}

\def \L {\mathcal{L}}

\def \O {\mathcal{O}}

\def \C {\mathbb{C}}

\def \P {\mathbb{P}}

\def \R {\mathbb{R}}

\def \Z {\mathbb{Z}}

\def \geq {\geqslant}

\def \kappa {\varkappa}

\def\={\;=\;}
\def\bal{\begin{aligned}}
\def\eal{\end{aligned}}

\DeclareMathOperator{\RHom}{RHom}

\providecommand{\arxiv}[1]{\href{http://arxiv.org/abs/#1}{arXiv:#1}}

\title{The conifold point}
\author{Sergey Galkin}
\address{National Research University Higher School of Economics, Faculty of Mathematics and Laboratory of Algebraic Geometry, Moscow}
\email{Sergey.Galkin@phystech.edu} 
\date{\today}
\begin{document}
\begin{abstract}
Consider a Laurent polynomial with real positive coefficients such that the origin is strictly inside its Newton polytope.
Then it is strongly convex as a function of real positive argument.
So it has a distinguished Morse critical point --- the unique critical point with real positive coordinates.

As a consequence we obtain a positive answer to a question of Ostrover and Tyomkin:
the quantum cohomology algebra of a toric Fano manifold contains a field as a direct summand.
Moreover, it gives a good evidence that the same statement holds for \emph{any} Fano manifold.

\end{abstract}
\maketitle
\setcounter{section}{0}


Let $W$ be a Laurent polynomial of variables $x_1,\dots,x_d$ with complex coefficients $a_n \in \C$:
\[ W = \sum_{n \in \Z^d} a_n x^n = \sum_{(n_1,\dots,n_d) \in \Z^d} a_{n_1,\dots,n_d} x_1^{n_1} \dots x_d^{n_d}. \]
Recall that the \emph{Newton polytope} $\Delta_W$ of a Laurent polynomial $W$ is defined as the convex hull in 
$\R^d$ of $n\in \Z^d$ such that $a_n \neq 0$. Assume further that $\Delta_W$ is $d$-dimensional
and the origin $0$ lies strictly inside $\Delta_W$.
In \cite{DK} Duistermaat and van der Kallen proved that there are infinitely many natural $k$ such that $k$-th moments $M_k(W)$
(defined as the constant terms of a Laurent polynomial $W^k$) does not vanish. Moreover, they proved that the generating function
$\hat{G}_W (t) = \sum_{k\geq 0} t^k M_k(W)$ has finite radius of convergence $0<R<\infty$, and function $\hat{G}_W (t)$
has logarithmic monodromy around some point $t_0$ with $|t_0| = R$. They have to use Hironaka's resolution of singularities in order to prove these theorems.
In this note we restrict to the case where coefficients are real and non-negative $a_n \in \R_+$. This condition greatly simplifies the picture,
existence of infinitely many non-vanishing moments is almost obvious, and additionally there is a very distinguished critical point $W$ that is non-degenerate.

Let $u_1,\dots,u_d$ be the coordinates on $\C^d$, the analytic map $x_i = \exp(u_i)$ is an \'etale $\Z^d$-covering $\exp: \C^d \to (\C^*)^d$. 
For a pullback
\[ \exp^* W = \sum_{n \in \Z^d} a_n e^{(u,n)} = \sum_{(n_1,\dots,n_d) \in \Z^d} a_{n_1,\dots,n_d} e^{(u_1 n_1 + \dots + u_d n_d)} \]
the partial derivation $\frac{\partial}{\partial u_i}$ coincides with the pullback of the logarithmic derivation $x_i \frac{\partial}{\partial x_i}$.
Note that 
$\exp$ establishes an isomorphism between
the domain 
$\R^d \subset \C^d$
(where all $u_i$ are real) 
and the domain 
$\R_+^d \subset (\C^*)^d$
(where all $x_i$ are real and positive),
in what follows we consider these identified domains as a topological manifold $T_+$.

\begin{lemma*}
Under the assumptions above:
\begin{enumerate}
\item 
The Hessian matrix $H_{ij} = \frac{\partial^2 W}{\partial u_i \partial u_j}$ is positive-semidefinite on $T_+$,
so $W : T_+ \to \R$ is a convex function.
\item Moreover 
$H_{ij}$ is positive-definite on $T_+$,
so $W: T_+ \to \R$ is strictly convex.
\item Moreover, 
the function $W: T_+ \to \R$ is strongly convex and attains a global minimum at some point $P \in T_+$.
\item Point $P$ is the unique critical point of $W$ in domain $T_+$, i.e. $dW_{u=u_0} = 0 \iff u_0 = P$ for any $u_0 \in T_+$.
\item Point $P$ is Morse i.e. the Hessian matrix of $W$ at $P$ is non-degenerate.
\end{enumerate}
In what follows we'll refer to $P$ as \emph{\bfseries the conifold point}.
\end{lemma*}
\begin{proof}

Clearly, partial derivatives of $W$ are given by
\[ \frac{\partial W}{\partial u_i} = \sum_n n_i a_n x^n \]
and
\[ H_{ij} = \frac{\partial^2 W}{\partial u_i \partial u_j} = \sum_n n_i n_j a_n x^n. \]

For any vector $v = (v_1,\dots,v_d)$ we have 
\[ v^t H(W) v = \sum_{i,j} v_i H_{ij} v_j = \sum_{i,j} \sum_n v_i (n_i n_j a_n x^n) v_j =
\sum_n a_n x^n (\sum_i v_i n_i) (\sum_j v_j n_j) = \sum_n a_n x^n (\sum_k v_k n_k)^2 \geq 0, \]
this proves the first statement.
Moreover, the last expression shows that $v^t H v = 0 \iff (v,n)=0$ for all $n$ such that $a_n \neq 0$.
So if vectors $n$ s.t. $a_n\neq0$ generate $\R^d$ then $v^t H v = 0 \iff v=0$, this proves the second statement.
Since $0$ is contained in the interior of the Newton polytope
for every direction $|u| \to \infty$ one of the monomials of $W$ also goes to $+\infty$.
Since all coefficients are positive $W$ goes to $+\infty$ as well.
This implies that $W$ has a global minimum $P$ on $T_+$,
similarly one proves that the minimal eigenvalue of $H(W)$ obtains a global minimum on $T_+$ so $W$ is strongly convex.
Point $P$ is critical because it is a minimum,
on the other hand strictly convex functions have at most one critical point.
Finally, the last statement (5) formally follows from (2): since Hessian matrix is positive-definite it is in particular non-degenerate.
\end{proof}

\begin{theorem*}[about toric Fano manifolds]
Let $X$ be a toric Fano manifold with a toric symplectic form $\omega$.
Then the small quantum cohomology algebra $QH(X,\omega)$ contains a field as a direct summand.
\end{theorem*}
\begin{proof}
By Proposition 3.3 of \cite{OT} the algebra $QH(X,\omega)$ coincides
with the Jacobi ring 
\[ J_W = \C[x_1^{\pm1},\dots,x_d^{\pm1}]/(\frac{d W}{d x_1},\dots,\frac{d W}{d x_d}) \]
of a particular combinatorially constructed Laurent polynomial $W_{X,\omega}$.
Moreover, as explained in Subsection 3.3 of loc.cit. to prove the Theorem for any $\omega$
it would suffice to consider the monotone case (i.e. $[\omega]=c_1(X)$). 
In the monotone case the respective Laurent polynomial equals 
$W_X = \sum_{v} x^v$,
here $v$ runs over all primitive generators of all rays of the fan of $X$,
and $x^v$ is the respective monomial. Clearly, $W_X$ satisfies the conditions of the Lemma,
so $W_X$ has a non-degenerate critical point $P$; point $P$ contributes a field as a direct summand of the Jacobi ring $J_W$.
\end{proof}

The Theorem above gives a positive answer to a question of Ostrover--Tyomkin \cite{OT}.
Their question was in turn raised as a modification of a question of Entov--Polterovich \cite{EP},
for which Proposition B of \cite{OT} gives a negative answer:
there are monotone symplectic toric Fano $4$-folds $(X,\omega)$ such that algebra $QH(X,\omega)$ is not semi-simple.
Nevertheless, for toric Fano manifolds $X$ the algebra $QH(X,\omega)$ is semi-simple for 
a \emph{generic} choice of $\omega$: see Corollary $5.12$ in \cite{I}, Proposition $7.6$ in \cite{FOOO} and Theorem A in \cite{OT}.
For non-toric Fano manifolds even generic semi-simplicity usually fails due to strong homological obstructions:
Theorem 1.8.1 of \cite{BM} and Theorem 1.3 of \cite{HMT} imply that if $QH(X,\omega)$ is semi-simple then $h^{p,q}(X) = 0$ unless $p=q$.
In contrast, next Remark aims to explain that the analogue of the Theorem should also hold for non-toric symplectic Fano manifolds.

{\bf Remark} regarding non-toric Fano manifolds.
The argument combines pictures of SYZ (Strominger-Yau-Zaslow \cite{SYZ}) and HMS (homological mirror symmetry \cite{K}) with ideas of \cite{H,CO,FOOO},
and very sketchy it goes as follows (see \cite{A,NNU} for at least some details).
Start from an arbitrary (non necessarily toric) Fano manifold $Y$.
Degenerate $Y$ to a (singular) toric Fano variety $X_0$.
The moment map $\mu: X_0 \to B$ gives a special Lagrangian tori fibration over the interior $\mu: X_0^o \to B^o$.
The symplectic transport (given by the distribution of the orthogonals of the fibers of degeneration with respect to symplectic form)
establishes a symplectomorphism between $X_0^o$ and an open subset in $Y$ thus producing many special Lagrangian tori $L_b \subset Y$.
The potential $W$ is then constructed as Fukaya-Oh-Ohta-Ono's obstruction $m_0(L_b,\nabla)$ for a (monotone) fiber $L_b$ equipped with a flat $U(1)$-connection $\nabla$,
that is a generating function for Maslov index two pseudoholomorphic discs bounded on $L_b$.
The case of a smooth toric Fano manifold $Y$ is explicitly computed in \cite{CO},
and the case of small degenerations in \cite{NNU}. The positivity of coefficients of $W$ is the geometrically evident fact,
that sometimes could be proved, e.g. if there are no pseudoholomorphic discs of negative Maslov index, or if those discs are away from $L_b$.
Thus by Lemma there is the conifold point $P$.
Section $6$ of \cite{A} and \cite{FOOO} explains how to identify the respective field summand of the Jacobi ring with a subalgebra of $QH(Y)$.
In particular, Theorem 6.1 of \cite{A} ensures that in smooth toric case, the set of all critical values of $W$
coincides with the set of eigenvalues of the quantum multiplication operator $\star_0 c_1(Y) : QH(Y)$.

The Theorem and the Remark above help to partially resolve the following conjecture of \cite{GGI}[(Section 3.1)].
For a Fano manifold $Y$ consider the set $U_Y$ of all eigenvalues $u_i$ of the quantum multiplication operator $\star_0 c_1(Y)$,
denote $T = \max |u_i|$. The {\bf Conjecture $\O$} says: $T$ lies in $U_Y$, for any $u_i \in U_Y$ if $|u_i| = T$ then $u_i/T$ is a root of unity,
and the multiplicity in $U_Y$ of the eigenvalue $T$ equals one.

We formulated it together with Hiroshi Iritani and Vasily Golyshev.
It appeared as a part of our investigation of the Gamma-conjectures (the relation between the (asymptotic) Apery class, the Gamma class,
and the integral structure in quantum cohomology), which turned out to be a close relative of Conjecture 4.2.2(3) of Dubrovin \cite{D}.

{\bf Number $T_Y$} is a real positive algebraic integer which can be considered as a symplectic invariant of a monotone Fano manifold $Y$.
The Lemma easily implies that if $Y$ is a toric Fano manifold, then $T_Y$ is bounded from above by $\dim Y + b_2(Y)$.
Non-toric Fano manifolds usually have $T_Y$ (high) above this bound. On the other hand, we were not able to prove any lower bound for $T_Y$,
even in the toric case. A plausible conjecture for the lower bound is $T_Y \geq \dim Y + 1$, with equality only for the projective space $Y \simeq \P^d$.

{\bf Relation to other work.}
The conifold point is explicitly constructed in the proof (due to Iritani) of Proposition 12.3 of \cite{NNU}.

Apparently it was observed by van Enckevort and van Straten as the Hypothesis 1 (H2) of \cite{ES}.
In some particular cases, such as for the mirror dual family of quintic threefolds
the terminology for the conifold point is well-established and its origin was not questioned.

The numerical conjecture $\O$ above constrains the geometry of the set $U_Y \subset \C$ of the \emph{critical values},
but does not say much about the set of the critical points. In case $W$ is mirror dual to a Fano manifold $Y$
we expect that $P$ is the unique singular point in the fiber $W^{-1} (W(P))$.
In contrast, such uniqueness sometimes fails for orbifolds, e.g. for the global quotient $\P^2/(\Z/3\Z)$
the mirror dual is $W = \frac{x_1^2}{x_2} + \frac{x_2^2}{x_1} + \frac{1}{x_1 x_2}$,
and it has exactly three conifold points over each of its three critical values.

In the discussion above we relate the $A$-model (quantum cohomology and other symplectic invariants) of the Fano manifold $Y$
to the $B$-model (Jacobi ring which can be thought as an algebro-geometric invariant) of the mirror-dual potential $W$.
Another direction of HMS relates the $B$-model of $Y$
(e.g. the bounded derived category of coherent sheaves $\D^b_{coh} Y$)
to the $A$-model of $W$ (the wrapped Fukaya category, or Fukaya--Seidel category of vanishing Lagrangian cycles).
Flip the direction of HMS to obtain the conjecture that explains the terminology:
the structure sheaf $\O_Y$ (as an exceptional object in $\D^b_{coh} Y$) corresponds under HMS to  
the real positive locus $T_+ \subset (\C^*)^d$ (considered as a Lefschetz thimble).
In this form we heard it from Denis Auroux.
Now folklore Conjecture $1$ in \cite{Ab}
is its relative, it can be deduced by combining SYZ, HMS, and the fact that the
equality $\RHom(\L,\O_y) = \C[0]$ for any structure sheaf $\O_y$ of a point $y\in Y$ implies that $\L$ is a line bundle.

{\bf Acknowledgement.}
We thank Ilya Tyomkin for explaining their motivations, the background and the history of their question.
The convexity argument is thanks to Grigory Mikhalkin.
We thank Denis Auroux, Hiroshi Iritani and Duco van Straten for interesting discussions.
Apparently all the ideas above are familiar to them.


\end{document}